\documentclass[12pt,reqno]{amsart}

\DeclareMathOperator{\ann}{ann}%
\DeclareMathOperator{\cd}{cd}%
\DeclareMathOperator{\cor}{cor}%
\DeclareMathOperator{\res}{res}%

\usepackage{amssymb,amscd}
\usepackage{hyperref}

\newcommand{\F}{\mathbb{F}}

\newcommand{\Fp}{\F_p}
\newcommand{\Gal}{\text{\rm Gal}}
\newcommand{\N}{\mathbb{N}}
\newcommand{\Z}{\mathbb{Z}}

\begin{document}

\title{Hilbert 90 for Galois Cohomology}

\begin{abstract}
    Assuming the Bloch-Kato Conjecture, we determine precise
    conditions under which Hilbert 90 is valid for Milnor $k$-theory
    and Galois cohomology.  In particular, Hilbert 90 holds for
    degree $n$ when the cohomological dimension of the Galois group
    of the maximal $p$-extension of $F$ is at most $n$.
\end{abstract}

\date{January 25, 2006}

\author[Lemire]{Nicole Lemire$^{\dag}$}
\address{Department of Mathematics, Middlesex College, \
University of Western Ontario, London, Ontario \ N6A 5B7 \ CANADA}
\thanks{$^\dag$Research supported in part by NSERC grant R3276A01.}
\email{nlemire@uwo.ca}

\author[Min\'{a}\v{c}]{J\'an Min\'a\v{c}$^{\ddag}$}
\thanks{$^\ddag$Research supported in part by NSERC grant R0370A01,
and by a Distinguished Research Professorship at the
University of Western Ontario.} \email{minac@uwo.ca}

\author[Schultz]{Andrew Schultz}
\address{Department of Mathematics, Building 380, Stanford
University, Stanford, California \ 94305-2125 \ USA}
\email{aschultz@stanford.edu}

\author[Swallow]{John Swallow}
\address{Department of Mathematics, Davidson College, Box 7046,
Davidson, North Carolina \ 28035-7046 \ USA}
\email{joswallow@davidson.edu}

\maketitle

\newtheorem{theorem}{Theorem}
\newtheorem{proposition}{Proposition}
\newtheorem{lemma}{Lemma}
\newtheorem{corollary}{Corollary}[theorem]
\newtheorem{corollaryc}[theorem]{Corollary}

\theoremstyle{definition}
\newtheorem*{remark*}{Remark}
\newtheorem{remark}{Remark}

\parskip=12pt plus 2pt minus 2pt

The key to the Bloch-Kato Conjecture is Hilbert 90 for Milnor
$K$-theory for cyclic extensions $E/F$ of degree $p$. It is
desirable to know when Hilbert 90 holds for Galois cohomology
$H^n(E,\Fp)$ as well.  In this paper we develop precise conditions
under which Hilbert 90 holds for Galois cohomology.

Let $p$ be a prime number, $E/F$ a cyclic extension of degree $p$
with group $G$, and assume that $F$ contains a primitive $p$th root
of unity $\xi_p$.  Let $\sigma$ be a generator of $G$. We say that
Galois cohomology satisfies Hilbert 90 for $n\in \N$ if the
following sequence is exact:
\begin{equation}\label{eq2}
    H^n(E,\Fp) \xrightarrow{\sigma-1} H^n(E,\Fp)
    \xrightarrow{\cor_{E/F}} H^n(F,\Fp).
\end{equation}
We say then that $(h90)_n$ is valid.

In section~\ref{se:conds} we formulate our results.  We recall some
basic results related to the Bloch-Kato Conjecture in
section~\ref{se:bkktheory}, and in sections~\ref{se:thms12} and
\ref{se:laterthms} we prove our results.  In section~\ref{se:hs} we
compare $(h90)_n$ with vanishing of the first cohomology of
$\Fp[G]$-module $H^n(E,\Fp)$.  We refer the reader to \cite{Lo} and
the references therein for interesting comments on classical Hilbert
90 and its extensions by Speiser and Noether.

\section{Conditions for Hilbert 90 for Galois
cohomology}\label{se:conds}

Choose $a\in F$ such that $E=F(\root{p}\of{a})$.  To simplify
notation we abbreviate $H^n(E,\Fp)$ and $H^n(F,\Fp)$ by $H^n(E)$ and
$H^n(F)$, respectively, and express cup products $\alpha\cup \beta$
as $\alpha.\beta$.

\begin{theorem}\label{th:p2}
    Suppose that $p=2$ and $n\in \N$.  Then $(h 90)_n$ is valid if
    and only if
    \begin{equation*}
        H^n(F) = \cor H^n(E) + (a).H^{n-1}(F).
    \end{equation*}
\end{theorem}

Let $E^\times:=E\setminus \{0\}$ and suppose still that $p=2$.  If
$a$ is a sum of two squares in $F$ then $a \in N_{E/F} E^\times$,
and by the projection formula (see \cite[Prop.~1.5.3iv]{NSW}) we
obtain $( a ). H^{n-1}(F) \subset \cor H^n(E)$. Therefore we have
the following corollary. (Observe that $a$ is a sum of two squares
in $F$ if and only if $-1 \in N_{E/F} E^\times$.)

\begin{corollary}\label{co:p2}
    Suppose that $p=2$, $n\in \N$, and $a\in F^2+F^2$. Then $(h
    90)_n$ is valid if and only if $\cor : H^n(E) \to
    H^n(F)$ is surjective.
\end{corollary}

\begin{theorem}\label{th:pnot2}
    Suppose that $p>2$ and $n \in \N$. Then $(h 90)_n$ is valid if
    and only if
    \begin{equation*}
        H^n(F) = \cor H^n(E) + ( \xi_p ).H^{n-1}(F).
    \end{equation*}
\end{theorem}

Now if $\xi_p \in N_{E/F} E^\times$, then by the projection formula,
$(\xi_p).H^{n-1}(F) \subset \cor H^n(E)$. Therefore we have
\begin{corollary}\label{co:pnot2}
    Suppose that $p>2$ and $\xi_p \in N_{E/F} E^\times$. Then $(h
    90)_n$ is valid if and only if $\cor : H^n(E) \to
    H^n(F)$ is surjective.
\end{corollary}

Corollaries~\ref{co:p2} and \ref{co:pnot2} together imply

\begin{corollaryc}\label{co:allp}
    Suppose that $n \in \N$ and $\xi_p \in N_{E/F} E^\times$. Then
    $(h90)_n$ is valid if and only if $\cor : H^n(E) \to
    H^n(F)$ is surjective.
\end{corollaryc}

Now let $F(p)$ be the compositum of all finite Galois $p$-power
extensions of $F$ in a fixed separable closure of $F$, and consider
the cohomological dimension of $G_{F}(p) := \Gal(F(p)/F)$. If $\cd
G_F(p) \le n$, then by \cite[Prop.~3.3.8]{NSW} we have the
surjectivity of $\cor:H^n(E) \to H^n(F)$. In fact, since $\xi_p\in
F$, the statement $\cd G_F(p)\le n$ is equivalent to the
surjectivity of $\cor$ on $H^n(E)$ for all cyclic extensions $E$ of
$F$ of degree $p$ \cite[Thm.~2]{LLMS}. Moreover, if a primitive
$p^2$th root of unity is in $F$, then for a suitable choice of the
$p^2$th root of unity $\xi_{p^2}\in F$ we obtain
$N_{E/F}(\xi_{p^2})=\xi_p$ for every cyclic extension $E/F$ of
degree $p$. Hence from Theorems~\ref{th:p2} and \ref{th:pnot2} and
Corollaries~\ref{co:p2} and \ref{co:pnot2} we have

\begin{corollaryc}\label{co:free}
    Suppose that $n\in \N$.  If $\cd G_F(p) \le n$, then $(h90)_n$
    is valid for all cyclic extensions $E/F$ of degree $p$.
    Moreover, if $\xi_{p^2}\in F$ and if $(h90)_n$ is valid for
    all cyclic extensions $E/F$ of degree $p$, then $\cd G_F(p) \le
    n$.
\end{corollaryc}

We present two further results: an interpretation of $(h90)_n$ in
terms of direct summands of $H^n(E)$ defined over $F$, and a
hereditary property of $(h90)_n$.

\begin{theorem}\label{th:basefield}
    Suppose that $n \in \N$. Then $(h 90)_n$ is valid if and only if
    for each pair $Q \subset H^n(F)$, $P \subset H^n(E)$ satisfying
    \begin{equation*}
        H^n(E) = \res Q \oplus P,
    \end{equation*}
    we have $\res Q = \{0\}$.
\end{theorem}

\begin{theorem}\label{th:hh90}
    Suppose that $n \in \N$.  Then if $(h 90)_n$ is valid, $(h
    90)_m$ is valid for all $m \ge n$.
\end{theorem}

At the end of the paper we compare $(h90)_n$ with
$H^1(G,H^n(E))=\{0\}$.

\section{Bloch-Kato and Milnor $K$-theory}\label{se:bkktheory}

For $i\ge 0$, let $K_iF$ denote the $i$th Milnor $K$-group of the
field $F$, with standard generators denoted by
$\{f_1\}\cdots\{f_i\}$, $f_1, \dots, f_i\in F^\times$. (See
\cite{Mi} and \cite[Chap.~IX]{FV}.) We use the usual abbreviation
$k_nF$ for $K_nF/pK_nF$.  For an extension of fields $E/F$, we use
$i_{E}$ for the natural inclusions of $K$-theory and $k$-theory, and
we also use $N_{E/F}$ for both the norm map $K_nE\to K_nF$ and the
induced map $k_nE\to k_nF$.

We prove our results first for Milnor $k$-theory, using Hilbert 90
for Milnor $K$-theory, and then we use the Bloch-Kato Conjecture to
identify $k_nF$ and $H^n(F)$. (See \cite[Lemma~6.11 and \S 7]{Vo1}
and \cite[\S 6 and Theorem~7.1]{Vo2}.  For further expositions of
the work of Rost and Voevodsky on Bloch-Kato Conjecture, see
\cite{Ro}, \cite{MVW}, and \cite{Su}.)

We say that Milnor $k$-theory satisfies Hilbert~90 at $n \in \N$ for
a cyclic extension $E/F$ of degree $p$ with $\Gal(E/F)=\langle
\sigma \rangle$ if the following sequence is exact:
\begin{equation}\label{eq1}
    k_n E \xrightarrow{\ \sigma-1\ } k_n E \xrightarrow{N_{E/F}} k_n
    F.
\end{equation}
We say then that $(h90)^M_n$ is valid. By the Bloch-Kato Conjecture,
there exists a $G$-equivariant isomorphism $k_nE\to H^n(E)$, and
therefore $(h90)^M_n$ is equivalent to $(h90)_n$. To determine
conditions for \eqref{eq2} it is then sufficient to determine
conditions for \eqref{eq1}.  Since the $G$-equivariant isomorphism
sends products to cup products, conditions for $k_nE$ expressed in
terms of products will carry over to the analogous conditions for
$H^n(E)$ expressed in terms of cup products.

We use the following two results in Voevodsky's work on the
Bloch-Kato Conjecture. Because we apply Voevodsky's results in the
case when the base field contains a primitive $p$th root of unity we
give formulations restricted to this case.  The first result is
Hilbert 90 for Milnor $K$-theory.

\begin{theorem}[{\cite[Lemma~6.11 and \S 7]{Vo1} and
\cite[\S 6 and Thm.~7.1]{Vo2}}] \label{th:bk}\
    Let $F$ be a field containing a primitive $p$th root of unity
    and $m\in\N$. For any cyclic extension $E/F$ of degree $p$
    with $\Gal(E/F)=\langle\sigma\rangle$, the
    sequence
    \begin{equation*}
        K_mE \xrightarrow{\ \sigma-1\ } K_mE \xrightarrow{N_{E/F}}
        K_mF
    \end{equation*}
    is exact.
\end{theorem}

As is standard, we then have the so-called ``Small Hilbert 90 for
$k_n$'':
\begin{corollary}\label{co:H90modp}
    Assume the hypotheses of Theorem~\ref{th:bk}.  Then the sequence
    \begin{equation*}
        (\sigma-1)k_mE+i_Ek_mF \longrightarrow k_mE
        \xrightarrow{N_{E/F}} k_mF
    \end{equation*}
    is exact.  In particular,
    \begin{equation*}
        (h90)_m^M \ \ \Longleftrightarrow
        \ \ i_E k_mF \subset (\sigma-1)k_m E.
    \end{equation*}
\end{corollary}

\begin{proof}
    Since $N_{E/F}(\sigma-1)= 0$ on $k_mE$ and $N_{E/F}$ is
    multiplication by $p$ on $i_Ek_nF$, our sequence is a complex.

    Now let $\alpha\in K_mE$ and write $\bar \alpha$ for the class
    of $\alpha$ in $k_mE$.  Suppose that $N_{E/F}\bar \alpha = 0\in
    k_mF$.  Then $N_{E/F}\alpha=p\beta$ for some $\beta\in K_mF$.
    Consider $\alpha'=\alpha-i_E(\beta)$. Then $N_{E/F}\alpha'=0$
    and by Theorem~\ref{th:bk} there exists $\gamma\in K_mF$ such that
    $(\sigma-1)\gamma=\alpha'=\alpha-i_E(\beta)$.  Then modulo $p$
    we have $(\sigma-1)\bar\gamma=\bar\alpha-i_E(\bar\beta)$ and so
    $\bar \alpha=(\sigma-1)\bar\gamma+i_E(\bar\beta)$.
\end{proof}

The following theorem is a strengthening of \cite[Prop.~5.2]{Vo1}.
Again $a\in F^\times$ is chosen so that $E=F(\root{p}\of{a})$.

\begin{theorem}[{\cite[Thm.~6]{LMS14}}]\label{th:esext}
    Let $F$ be a field containing a primitive $p$th root of unity.
    Then for any cyclic extension $E/F$ of degree $p$ and $m\in \N$
    the sequence
    \begin{equation*}
        k_{m-1}E \xrightarrow{N_{E/F}} k_{m-1}F \xrightarrow{\{a\}
        \cdot -\ } k_m F \xrightarrow{i_E} k_m E
    \end{equation*}
    is exact.
\end{theorem}

\section{Proofs of Theorems~\ref{th:p2} and
\ref{th:pnot2}}\label{se:thms12}

\begin{proof}[Proof of Theorem~\ref{th:p2}]
    By Corollary~\ref{co:H90modp}, $(h90)_n^M$ is valid if and only
    if
    \begin{equation*}
        i_E k_nF\subset (\sigma-1) k_nE = (\sigma+1) k_n E = i_E
        N_{E/F} k_n E.
    \end{equation*}
    Hence $(h90)_n^M$ is valid if and only if $k_n F = N_{E/F}k_n E
    + \ker i_E$.  By Theorem~\ref{th:esext}, $\ker i_E = \{a\}\cdot
    k_{n-1}F$.  Therefore
    \begin{equation*}
        (h90)^M_n \ \ \Longleftrightarrow\ \ k_n F = N_{E/F}k_n E +
        \{a\}\cdot k_{n-1}F.
    \end{equation*}
    Using the Bloch-Kato Conjecture, the theorem follows.
\end{proof}

\begin{remark}
    The theorem also follows from Arason's long exact sequence (see
    \cite{A}).  Therefore when $p=2$ we can deduce $(h90)_n$ from
    basic Galois cohomology, without reference to the Bloch-Kato
    Conjecture.
\end{remark}

\begin{proposition}~\label{pr:sigmamin1}
    Let $F$ be a field containing a primitive $p$th root $\xi_p$ of
    unity and $n\in\N$. For any cyclic extension $E/F$ of degree $p$
    with $\Gal(E/F)=\langle\sigma\rangle$,
    \begin{equation*}
        (\sigma-1)k_nE\cap (k_nE)^G=i_{E}(\{\xi_p\}\cdot k_{n-1}F) +
        i_{E}N_{E/F}k_nE.
    \end{equation*}
\end{proposition}

Before proving Proposition~\ref{pr:sigmamin1} we introduce some
further notation and establish Lemma~\ref{le:keycharnotp} below.

For $y\in k_nE\setminus \{0\}$, set the length $l(y)$ of $y$ to be
\begin{equation*}
    l(y):=\max\{t\in \N : (\sigma-1)^{t-1} y \neq 0\}.
\end{equation*}
Observe that since $1+\sigma+\cdots+\sigma^{p-1}= (\sigma-1)^{p-1}$
in $\Fp G$, $i_E N_{E/F}= (\sigma-1)^{p-1}$ on $k_nE$ and therefore
$l(y)\le p$.

Since $\{\xi_p\}\cdot k_{n-1}F = \{\xi_p^c\}\cdot k_{n-1}F$ for
$(p,c)=1$, we assume without loss of generality for the proofs of
the proposition and the following lemma that $\root{p}
\of{a}^{\sigma-1} = \xi_p$.

\begin{lemma}\label{le:keycharnotp}
    Suppose $y\in k_nE$ with $l=l(y)\ge 2$. Then if $l\ge 3$,
    \begin{equation*}
        (\sigma-1)^{l-1}(y)\in i_{E}N_{E/F}k_nE,
    \end{equation*}
    and if $l=2$
    \begin{equation*}
        (\sigma-1)y\in i_{E}(\{\xi_p\}\cdot k_{n-1}F+ N_{E/F}k_nE).
    \end{equation*}
\end{lemma}

\begin{proof}
    If $l=p$, then $(\sigma-1)^{p-1}k_nE=i_{E}N_{E/F}k_nE$ shows the
    result in this case.  Hence we may assume that $p>2$.

    Suppose $l<p$. Then $y\in \ker(\sigma-1)^{p-1}$ and so
    $i_{E}N_{E/F}(y)=0$.  By Theorem~\ref{th:esext}, there exists
    $b\in k_{n-1}F$ such that $N_{E/F}y=\{a\}\cdot b$. By the
    projection formula~\cite[p.~81]{FW},
    \begin{equation*}
        N_{E/F}\big(y-\{\root{p}\of{a}\}\cdot i_{E}(b)\big)=0.
    \end{equation*}
    Then by Corollary~\ref{co:H90modp}, there exist $\omega\in
    k_nE$ and $f\in k_nF$ such that
    \begin{equation*}
        y=(\sigma-1)\omega+\{\root{p}\of{a}\}\cdot i_{E}(b)+
        i_{E}f
    \end{equation*}
    and hence
    \begin{equation*}
        (\sigma-1)^{l-1}y=(\sigma-1)^l\omega+(\sigma-1)^{l-2}i_{E}
        (\{\xi_p\}\cdot b).
    \end{equation*}

    If $l(y)\ge 3$ we deduce
    \begin{equation*}
        (\sigma-1)^{l-1}y=(\sigma-1)^l\omega
    \end{equation*}
    where $l(\omega)=l+1$. Set $y_{l+1}=\omega$ and repeat the
    argument. We obtain $y_k\in k_nE$ of lengths $l\le k\le p$ with
    \begin{equation*}
        (\sigma-1)^{l-1}y=(\sigma-1)^{k-1}y_k.
    \end{equation*}
    Take $\alpha=y_p$ to obtain
    \begin{equation*}
        (\sigma-1)^{l-1}y\in i_{E}N_{E/F}k_nE,
    \end{equation*}
    as required.

    If $l(y)=2$ we have that
    \begin{equation*}
        (\sigma-1)y=(\sigma-1)^2\omega+i_{E}(\{\xi_p\}\cdot b)
    \end{equation*}
    for some $\omega\in k_nE$ and some $b\in k_{n-1}F$. We see
    that $l(\omega)\le 3$. If $l(\omega)<3$ then
    \begin{equation*}
        (\sigma-1)y\in i_{E}(\{\xi_p\}\cdot k_{n-1}F),
    \end{equation*}
    while if $l(\omega)=3$ then by the previous case we see that
    $(\sigma-1)^2 \omega\in i_{E}N_{E/F}k_nE$, so the result
    holds in either case.
\end{proof}

\begin{proof}[Proof of Proposition~\ref{pr:sigmamin1}]
    The right-hand side is contained in $i_{E}k_nF$ and so is fixed
    by $G$. Furthermore,
    \begin{equation*}
        i_{E}N_{E/F}k_nE = (\sigma-1)^{p-1}k_nE\subset
        (\sigma-1)k_nE.
    \end{equation*}
    Let $f\in k_{n-1}F$. Then because
    \begin{equation*}
        (\sigma-1)(\{\root{p}\of{a}\}\cdot
        i_{E}(f))=i_{E}(\{\xi_p\}\cdot f),
    \end{equation*}
    we have $i_{E}(\{\xi_p\}\cdot k_{n-1}F)\subset (\sigma-1)k_nE$,
    and so
    \begin{equation*}
        i_{E}(\{\xi_p\}\cdot k_{n-1}F)+i_{E}N_{E/F}k_nE\subset
        (\sigma-1)k_nE\cap (k_nE)^G.
    \end{equation*}
    Now let $0\ne x\in (k_nE)^G\cap (\sigma-1)k_nE$. Then
    $x=(\sigma-1)\gamma$ for some $\gamma\in k_nE$ with
    $l(\gamma)=2$. The result follows from
    Lemma~\ref{le:keycharnotp}.
\end{proof}

\begin{proof}[Proof of Theorem~\ref{th:pnot2}]
    By Corollary~\ref{co:H90modp}, $(h90)^M_n$ is valid if and only
    if $i_Ek_nF\subset (\sigma-1)k_nE$, or, equivalently,
    \begin{equation*}
        (h90)_n^M \ \ \Longleftrightarrow \ \ (\sigma-1)k_nE\cap
        i_Ek_nF=i_Ek_nF.
    \end{equation*}
    By Proposition~\ref{pr:sigmamin1},
    \begin{equation*}
        (\sigma-1)k_nE\cap i_Ek_nF=i_E(\{\xi_p\}\cdot
        k_{n-1}F)+i_EN_{E/F}k_nE.
    \end{equation*}
    Hence $(h90)^M_n$ is valid if and only if
    \begin{equation*}
        i_Ek_nF=i_E(\{\xi_p\}\cdot k_{n-1}F+N_{E/F}k_nE).
    \end{equation*}
    Since $p>2$, $\{a\}\cdot \{a\}=0$. Then by
    Theorem~\ref{th:esext}, $\ker i_E \subset N_{E/F}k_nE$. Hence
    \begin{equation*}
        (h90)_n^M \ \ \Longleftrightarrow \ \ k_nF = \{\xi_p\}\cdot
        k_{n-1}F+N_{E/F}k_nE.
    \end{equation*}
    Using the Bloch-Kato Conjecture, the theorem follows.
\end{proof}

\section{Proofs of Theorems~\ref{th:basefield} and
\ref{th:hh90}}\label{se:laterthms}

We recall some results on $\Fp[G]$-modules for $G$ a cyclic group of
order $p$.  The indecomposable $\Fp[G]$-modules are precisely the
cyclic $\Fp[G]$-modules $V_i:=\Fp[G]/\langle(\sigma-1)^{i}\rangle$
of dimensions $i=1,\dots, p$; hence $V_i$ is annihilated by
$(\sigma-1)^i$ but not by $(\sigma-1)^{i-1}$.  The trivial
$\Fp[G]$-module $V_1\simeq \Fp$ is the unique simple $\Fp[G]$-module
up to isomorphism.  Recall that a semisimple module is any direct
sum, possibly infinite, of simple modules. For each $i$, the fixed
submodule $V_i^G$ of $V_i$ is $(\sigma-1)^{i-1}V$, and each $V_i$
has finite composition length and therefore (see
\cite[Lemma~12.8]{AnFu}) its endomorphism ring is local.

\begin{proposition}\label{pr:maxtriv}
    Let $M$ be an $\Fp[G]$-module.  Then $M=S\oplus T$ where
    $S$ is a maximal semisimple direct summand and $T^G\subset
    (\sigma-1)T$.
\end{proposition}

\begin{proof}
    Since $\Fp[G]$ is an Artinian principal ideal ring, every
    $\Fp[G]$-module $W$ is a direct sum of cyclic $\Fp[G]$-modules
    \cite[Thm.~6.7]{SV}. By the Krull-Schmidt-Azumaya Theorem (see
    \cite[Thm.~12.6]{AnFu}), all decompositions of $W$ into direct
    sums of indecomposable modules are equivalent.  Now decompose
    $M=\oplus M_i$ where for each $i=1, \dots, p$, $M_i$ is a direct
    sum of cyclic $\Fp[G]$-modules $V_i$ of dimension $i$. Set $S=M_1$,
    $T=M_2\oplus \cdots\oplus M_p$.  Then $S$ is a semisimple direct
    summand and $T$ contains no nonzero semisimple direct summand.
    Hence $S$ is maximal.  Since for each $i\ge 2$, $T_i$ is a
    direct sum of modules $V_i$ satisfying $V_i^G =
    (\sigma-1)^{i-1}V_i$, we have that $T$ satisfies $T^G\subset
    (\sigma-1)T$.
\end{proof}

\begin{proof}[Proof of Theorem~\ref{th:basefield}]
    Assume first that $(h90)^M_n$ is valid, and suppose that $k_nE =
    i_E Q \oplus P$ for $Q\subset k_nF$, $P\subset k_nE$. Then
    $N_{E/F}i_E Q = pQ = \{0\}$.  By $(h90)^M_n$,
    \begin{equation*}
        i_E Q\subset (\sigma-1)k_nE = (\sigma-1)i_E Q \oplus
        (\sigma-1)P = (\sigma-1)P \subset P.
    \end{equation*}
    Hence $i_E Q \subset i_E Q\cap P = \{0\}$, as desired.

    For the other direction, assume that if $k_nE=i_E Q \oplus P$
    for $Q\subset k_nF$, $P\subset k_nE$, then $i_E Q = \{0\}$. By
    Proposition~\ref{pr:maxtriv} we may write $k_nE=S\oplus T$ where
    $S$ is a maximal semisimple direct summand and $T^G\subset
    (\sigma-1)T$.  We claim that $U=i_E k_n F\cap S=\{0\}$, as
    follows.  Because $S$ is semisimple there exists $V$ such that
    $U\oplus V=S$.  Let $Q\subset k_nF$ be the inverse image of $U$
    under $i_E$. Then $k_nE=i_EQ\oplus (V\oplus T)$ so that
    $U=\{0\}$, as desired.

    Now let $f\in k_nF$ be arbitrary, and write $x := i_E f = s+t$
    along $S\oplus T$.  Since $(\sigma-1)x=0$, we see that
    $(\sigma-1)t=0$ and $t\in T^G\subset (\sigma-1)T$.  By
    Proposition~\ref{pr:sigmamin1}, $t\in i_E k_n F$ and therefore
    $s=x-t\in i_E k_n F$.  But since $U=\{0\}$, $s=0$.  Hence $i_E
    k_n F\subset (\sigma-1)k_nE$.  By Corollary~\ref{co:H90modp}, we
    deduce that $(h90)^M_n$ is valid.

    Using the Bloch-Kato Conjecture, the theorem follows.
\end{proof}

\begin{remark}
    The same argument as in the first paragraph of the proof shows
    directly that if $(h90)_n$ is valid and $H^n(E) = \res Q \oplus
    P$ for $Q\subset H^n(F)$ and $P\subset H^n(E)$, then $\res Q =
    \{0\}$.
\end{remark}

\begin{proof}[Proof of Theorem~\ref{th:hh90}]
    Let $n \in \N$ and assume first that $p>2$. We prove the result
    by induction on $m$.  The base case $m=n$ is given. Assume then
    that $(h90)^M_m$ holds.  By Theorem~\ref{th:pnot2},
    \begin{equation*}
        k_m F = N_{E/F} k_m E + \{ \xi_p \} \cdot k_{m-1} F.
    \end{equation*}
    Then each element in $k_{m+1} F$ is a sum of the form $\{ f \}
    \cdot r + \{ g \} \cdot \{ \xi_p \} \cdot s$ where $f, g \in
    F^\times$, $r \in N_{E/F} k_m E$, and $s \in k_{m-1} F$.

    By the projection formula (see \cite[Thm.~3.8]{FV}), $\{ f \}
    \cdot r \in N_{E/F} k_{m+1} E$.  Moreover, $\{ g \} \cdot \{
    \xi_p \} \cdot s = \{ \xi_p \} \cdot \{ g^{-1} \} \cdot s \in \{
    \xi_p \} \cdot k_m F$. Therefore
    \begin{equation*}
        k_{m+1} F = N_{E/F} k_{m+1} E + \{ \xi_p \} \cdot k_m F.
    \end{equation*}
    By the proof of Theorem~\ref{th:pnot2}, we have $(h90)^M_{m+1}$.
    Using the Bloch-Kato Conjecture, the theorem for $p>2$ follows.

    The case $p=2$ follows by replacing $\xi_p$ with $a$ and
    Theorem~\ref{th:pnot2} with Theorem~\ref{th:p2} in the argument
    above.
\end{proof}

\section{Hilbert versus Noether and Speiser}\label{se:hs}

It is interesting to compare $(h 90)_n$ with the condition
$H^1(G,H^n(E)) = \{ 0 \}$.  Write $\ann_n (a)$ for the annihilator
of the cup-product with $(a)$ in $H^n(F)$.

\begin{theorem}\label{th:hs}
    Suppose $n\in \N$.  Then
    \begin{equation*}
        H^1(G,H^n(E))=\{0\}\ \ \Longrightarrow\ \ (h90)_n.
    \end{equation*}

    If $p>2$ and $n\in \N$, the following are equivalent:
    \begin{enumerate}
        \item $H^1(G,H^n(E))=\{0\}$
        \item $H^n(E)$ is a free $\Fp[G]$-module
        \item $\cor:H^{n-1}(E)\to H^{n-1}(F)$ is surjective.
    \end{enumerate}

    If $p=2$ and $n\in \N$, the following are equivalent:
    \begin{enumerate}
        \item $H^1(G,H^n(E))=\{0\}$
        \item $H^n(E)$ is a free $\Fp[G]$-module
        \item $H^n(E) = \ann_n (a) \oplus (a).H^{n-1}(F)$.
    \end{enumerate}
\end{theorem}

\begin{proof}
    If $H^1(G,H^n(E))=\{0\}$ then $(\sigma-1)H^n(E)=\ker N$, where
    the endomorphism $N:H^n(E)\to H^n(E)$ is defined by
    $N=\res \cor$.  Since $\ker N\supseteq \ker \cor$, we obtain
    $(h90)_n$.

    The equivalences (1)$\Leftrightarrow$(2) follow from \cite[\S
    III.1, Prop.~1.4]{L}.  The equivalences (2)$\Leftrightarrow$(3)
    follow from \cite[Thm.~1]{LMS2}.
\end{proof}

\section{Acknowledgements}

Andrew Schultz would like to thank Ravi Vakil for his encouragement
and direction in this and all other projects.  John Swallow would
like to thank Universit\'e Bordeaux I for its hospitality during
2005--2006.

\end{document}